\documentclass[]{elsarticle}
\usepackage{amsfonts,amsmath,amssymb}
\usepackage{multirow}
\usepackage{graphicx}
\usepackage{algorithmicx}
\usepackage[linesnumbered,ruled]{algorithm2e}
\usepackage{lmodern}
\usepackage{color}
\usepackage{makecell}
\usepackage{soul}

\newtheorem{theorem}{Theorem}

\newcommand{\R}{{\mathbb R}}

\newcommand{\K}{{\mathbb K}}

\newcommand{\bx}{\mbox{\boldmath{$x$}}}
\newcommand{\bb}{\mbox{\boldmath{$b$}}}
\newcommand{\bc}{\mbox{\boldmath{$c$}}}

\newcommand{\be}{\mbox{\boldmath{$e$}}}

\newcommand{\edit}[1]{\textcolor{black}{#1}}

\begin{document}
\begin{frontmatter}
\title{Approximating matrix functions by block\\ Krylov methods with randomized vectors\\
~\\
\small{Dedicated to Paul Van Dooren on the occasion of his 75th birthday.}   \\
}
	
\author{Josh Kane\fnref{ksu}}
\ead{jkane21@kent.edu}
\author{Lucas Onisk\fnref{emory}}
\ead{lonisk@emory.edu}
\author{Lothar Reichel\fnref{ksu}}
\ead{reichel@math.kent.edu}
\author{Giuseppe Rodriguez\fnref{unica}}
\ead{rodriguez@unica.it}

\address[ksu]{Department of Mathematical Sciences, Kent State University, 1300 Lefton 
Esplanade, Kent, OH 44242, USA}
\address[emory]{Department of Mathematics, Emory University, 400 Dowman Drive, Atlanta,
GA 30322, USA}
\address[unica]{Department of Mathematics and Computer Science, University of Cagliari, 
via Ospedale 72, Cagliari, 09124, Italy}

\begin{abstract}
The need to evaluate expressions of the form $f(A)\bb$, where $A$ is a square matrix, $f$
is a function, and $\bb$ is a vector, arises in several areas of applied mathematics. When
the matrix $A$ is very large, it is usually not attractive to evaluate $f(A)$. Instead, 
$f(A)\bb$ often is approximated by computing an estimate in a Krylov subspace that depends
on $A$ and $\bb$, and only requires that $f$ be evaluated at a small matrix. This paper 
explores the application of several variants of randomized block Krylov methods to the approximation 
of $f(A)\bb$. \edit{Computed examples suggest that block Krylov methods with an 
initial block vector that contains $\bb$ as well as a few randomly generated 
vectors may require less computing time and reduce the number of Krylov 
steps than standard Krylov methods.}
\end{abstract}

\end{frontmatter}

\graphicspath{{./figs/}}

\section{Introduction} \label{sec:1}
Many applications in scientific computing require the evaluation of expressions of the
form
\begin{equation}\label{fab}
f(A)\bb,
\end{equation}
where $A\in{\R}^{n\times n}$ is a large matrix, $f$ is a function such that $f(A)$ is well
defined, and $\bb\in\R^{n}$ is a vector. Applications include the solution of linear 
discrete ill-posed problems~\cite{CLR,CR02,CRZ,CDD}, network
analysis~\cite{BB20,DMR,EH10}, and the solution of linear partial differential
equations~\cite{DH05,GS,La}. It can be prohibitively expensive to compute
\eqref{fab} by 
first evaluating $f(A)$ when the matrix $A$ is large. Therefore, the use of Krylov 
subspace methods for the approximation of expressions of the form \eqref{fab} have 
received considerable attention. These methods first reduce the matrix $A$ to a small 
matrix $\mathcal{H}_{m,m}\in\R^{m\times m}$ with $1\leq m\ll n$ by computing a Krylov subspace basis of 
dimension $m$, which involves \edit{the evaluation of} $m$ matrix-vector products with $A$, and then compute
$f\left(\mathcal{H}_{m,m}\right)$. A thorough treatment of methods for the evaluation of $f$ of a small matrix is
provided by Higham~\cite{H08}; see also~\cite{BR,DK,HL} for further discussions.

Arithmetic floating point operations (FLOPs) are often used as a measure of
complexity of an algorithm. However, on many available computers, the FLOP count is not proportional to the required computing time. This holds
true not only for \edit{large} parallel computers, but also for personal 
\edit{laptop} computers with modest computing 
capabilities. For instance, in many computing environments, the simultaneous evaluation of
the product of a matrix with a few vectors is only insignificantly more time-consuming 
than the evaluation of the product of the same matrix with only one vector. This 
observation suggests that it may be beneficial to apply block Krylov methods instead of
standard Krylov methods (with block size one) to reduce the matrix $A$ to a small matrix, 
that is subsequently used to compute an approximation of \eqref{fab}. In fact, some 
authors propose to count the number of block vector products carried out by an algorithm
as its measure of complexity, independently of the block size used; see, e.g., Chen et al.~\cite{CHLZ}.

The use of randomized block Krylov methods for the determination of a few of the largest 
eigenvalues and associated eigenvectors of a large symmetric matrix, and for the 
computation of a few of the largest singular values and associated singular vectors of a 
large nonsymmetric, possibly rectangular, matrix have recently been discussed by 
Tropp and Webber~\cite{Tr,TW}. A discussion on the application of randomized block
Krylov methods to the Tikhonov problem can be found in~\cite{KOR}. This paper 
considers the application of a suitable number of steps of a randomized block 
Arnoldi or block Lanczos methods to reduce the large square matrix $A$ in 
\eqref{fab} to a small square matrix, which is used to determine an 
approximation of \eqref{fab}. We note that recently Persson et al.~\cite{Chen25} discussed the application of 
randomized block Krylov methods to the approximation of matrix functions $f(A)$ when $A$ 
is a symmetric matrix. This problem differs from the approximation of \eqref{fab} in that
the latter expression only requires the evaluation of a vector, which
depends on $\bb$. We will see that the presence of the vector $\bb$ in
\eqref{fab} is important for the design of randomized block Krylov methods.

Several approaches using randomized linear algebra have been used to speed up
the evaluation of expressions of the form \eqref{fab} in recent literature. For
instance, Cortinovis et al.~\cite{CKN24} sought to speed up the evaluation of
these expressions by dispensing with the
requirement that the basis vectors of the Krylov subspace be orthonormal. Specifically, randomization is introduced
by applying sketching to compute simplified approximations of the inner products that have
to be evaluated in standard Krylov methods (with block size one).
Guidotti et al.~\cite{GMAM25} also employ sketching to speed up the computations and so do G\"uttel and 
Schweitzer~\cite{GS22}. The methods of the present paper can be combined with sketching. 
Preliminary numerical experiments suggest that sketching may yield speed-ups for certain
expressions \eqref{fab}, but not for others. In this work, we focus on the performance of
randomized block Krylov methods that do not employ sketching.

Non-randomized block Arnoldi methods for matrix function approximation have fairly 
recently been studied by Frommer et al.~\cite{FLS}, who consider a framework that 
incorporates various inner products and scaling functions for block vector products (see Section \ref{sec:unified} for more detail). They then define the B(FOM)$^2$ method, which leverages this framework to approximate matrix functions. This framework includes classical, global, and loop-interchange block strategies. Global block methods were first discussed by Jbilou et al.~\cite{Jbilou99,JST}. A comparison of
global and classical block methods in the context of solving large linear discrete 
ill-posed problems that arise in color image restoration is furnished by Bentbib et al.~\cite{BEGJOR}. 

Our focus is to investigate the performance of 
alternate versions of the classical, global, and loop-interchange block methods that use randomized initial vectors to approximate expressions of the form 
\eqref{fab}.
Up to now, we have discussed
the approximation of \eqref{fab} with a single vector 
$\bb$. However, it is straightforward to replace this vector with a block 
vector with a few columns \edit{whose resulting algorithms we will refer to throughout as randomized block Krylov methods}.  Our numerical results show that using randomized block Krylov methods with block size larger than one, but not too large, may reduce both the number of steps 
required to obtain a specified tolerance by the block Krylov method and the total runtime. 

This paper is organized as follows: Section~\ref{sec:2} describes randomized 
block Krylov algorithms as well as the inner products and scalings to be used. 
Section~\ref{sec:3} discusses the approximation of expressions \eqref{fab} by these methods. 
Computed examples are presented in Section~\ref{sec:4} and Section~\ref{sec:5} contains 
concluding remarks.

Much work in numerical linear algebra is devoted to exploiting structure of computational
problems to reduce the computational burden or the sensitivity to errors in the data
to round-off errors introduced during computations. Paul Van Dooren has made many
significant contributions in this area including \cite{BMVD,BGHSVD,LMVD,Mas1}. The 
present paper exploits the structure in the matrix function approximation problem 
considered by applying block Krylov methods with the inclusion of randomized vectors. 

\section{Randomized block Krylov subspace methods}\label{sec:2}
\edit{In this section we begin with a review of the classical block Arnoldi and block Lanczos subspace methods that use an initial randomized block. These methods are applied by Tropp and Webber~\cite{Tr,TW} to compute low-rank approximations of large square matrices $A\in\R^{n\times n}$. Our approach in this work focuses on the modification and size of the initial block when approximating expressions of the form \eqref{fab}. Alternative block strategies that have shown competitive advantages in other settings compared to the classical strategy, namely, the global and loop-interchange strategies discussed by Frommer et al. in \cite{FLS} are then summarized in Section~\ref{global_LI}. We conclude in Section \ref{sec:unified}  with a unifying framework for all three block strategies considered.}

\subsection{Randomized classical block Arnoldi and Lanczos methods}\label{RBA}
We first consider the situation when the large matrix $A\in{\R}^{n\times n}$ in 
\eqref{fab} is nonsymmetric, and outline the application of the randomized classical 
block Arnoldi (RCBA) method with block size $1<\ell\ll n$. Then we describe the formation of the initial block \edit{that contains a multiple of the vector 
$\bb$ and $\ell-1$ random columns}. We conclude with a brief comment about the 
corresponding randomized classical block Lanczos (RCBL) method, which is used when the 
matrix $A$ is symmetric.

Let $V_1\in{\R}^{n\times \ell}$ be the initial block vector for the RCBA 
method.
The following discussion, as well as Algorithm~\ref{RCBA}, also applies when 
$\ell=1$, but in this case the algorithm can be simplified. The columns of $V_1$ are 
chosen to be orthonormal. Specifically, first let the entries of the matrix 
$\Omega_\ell\in\R^{n\times\ell}$ be normally distributed random numbers with mean zero and
variance one. Compute its ``skinny'' QR factorization $\Omega_\ell=Q_\ell R_{\ell}$, where
the matrix $Q_\ell\in\R^{n\times\ell}$ has orthonormal columns and 
$R_\ell\in\R^{\ell\times\ell}$ is upper triangular. We tacitly assume that such a QR factorization exists and that $R$ is nonsingular; this is the generic situation. We obtain a randomized
block Arnoldi method by letting 
\begin{equation}\label{startBlock0}
V_1=Q_\ell.
\end{equation}
Performing $q\ll n/\ell$ steps (which we refer to as the \emph{depth}) of the
RCBA algorithm applied to $A$ with initial
block vector $V_1$, as described in Algorithm~\ref{RCBA}, produces the decomposition
\begin{equation} \label{RBArelation}
        A\mathcal{V}_q = \mathcal{V}_{q+1}\mathcal{H}_{q},
\end{equation}
where the matrix $\mathcal{V}_{q+1} = \left[V_1\,V_2,\, \dots\,V_{q+1}\right] \in 
{\R}^{n\times \ell(q+1)}$ has orthonormal block columns $V_i\in{\R}^{n\times\ell}$ as well
as orthonormal columns. We consider two block columns $V_i$ and $V_j$ 
orthonormal if
\[
V_i^TV_j=\begin{cases} I_\ell, & \text{for } i= j,\\
		  O_\ell, & \text{for } i\ne j,
                \end{cases}
\]
where $I_\ell$ and $O_\ell$ denote the identity and zero matrices of order $\ell$,
respectively.

The matrix $\mathcal{V}_q\in{\R}^{n\times \ell q}$ is made up of the first $q$ 
block columns $V_i$ of $\mathcal{V}_{q+1}$ and the matrix
\begin{equation*}
\mathcal{H}_{q} = \begin{bmatrix}
    & H_{1,1}  & H_{1,2} &       &\ldots            &H_{1,q} \\
    & H_{2,1}  & H_{2,2} &       &           &  \\
    & 0              & H_{3,2} &       &\ddots            & \vdots  \\
    & \vdots         &\ddots        & &\ddots            &  \\
    &              &        &      & H_{q,q-1}  & H_{q,q}  \\
    & 0 & \cdots & \cdots & 0 & H_{q+1,\,q}
\end{bmatrix} \in \mathbb{R}^{\ell(q+1) \times \ell q}
\end{equation*}
is of upper block Hessenberg form. The nontrivial blocks $H_{i,j}\in\R^{\ell\times\ell}$
are defined in Algorithm~\ref{RCBA}. The subdiagonal blocks $H_{j+1,j}$ for $j=1,2,\dots,q$ are upper 
triangular matrices determined by QR factorizations. We discuss this further in Section~\ref{sec:unified}. \edit{The classical block Arnoldi process may be used to compute an approximation of the solution of linear systems of equations, 
whose right-hand side is made up of a block vector with $\ell>1$ columns. In 
this paper, we apply block Krylov subspace methods to approximate the 
expression \eqref{fab}. We refer to Bentbib et al. \cite{BEGJOR} and Saad~\cite[Chapter 6]{Sa1} for more details on the classical block Arnoldi process. }

\begin{algorithm}[ht!] \label{RCBA}
\caption{RCBA method}
\textbf{Input:} $A\in\R^{n\times n}$, $\bb\in\R^{n}$, block size $\ell \in \mathbb{N}$, number of steps $q \in \mathbb{N}$\;
\textbf{Output:} $\mathcal{H}_{q}\in\R^{\ell (q+1)\times \ell q}$, $\mathcal{V}_{q+1}\in\R^{n\times(q+1)\ell}$ which satisfy \eqref{RBArelation}\;
    Generate \edit{$\Omega_{\ell-1} \in\R^{n\times (\ell-1)}$ with $\mathcal{N}\left(0,1\right)$ entries}\;
    Prepend $\bb$ to $\edit{\Omega_{\ell-1}}$:\quad\edit{$\Omega_{\ell}= \left[\bb,\,\Omega_{\ell-1}\right]$}\;
    Normalize: set $[V_1,\sim]=\operatorname{qr}(\edit{\Omega_{\ell}},\textrm{`econ'})$ \edit{or via MGS}\;
		\For{$k=1,\ldots,q$}{
			Compute $W=AV_k$\;
				\For{$j=1,\ldots,k$}{
					$H_{j,k}=V_j^TW$\;
					$W=W-V_jH_{j,k}$\;\label{line11}
				}
			$[V_{k+1},H_{k+1,k}]=\operatorname{qr}(W,\textrm{`econ'})$\;
		}
\end{algorithm}

We have found the inclusion of the vector $\bb$ in the \edit{span of the} initial block vector 
$V_1$ to be beneficial when approximating expressions of 
the form \eqref{fab}. We therefore define an alternative initial block vector $V_1$ in 
Algorithm~\ref{RCBA} as
\begin{equation} \label{startBlock}
V_1 = \left[\bb,\,\,\widehat{Q}_{\ell-1} \right],
\end{equation}
where the matrix $\widehat{Q}_{\ell-1}\in\R^{n\times(\ell-1)}$ is determined as follows. 
Let $\Omega_{\ell-1}\in{\R}^{n\times(\ell-1)}$ be a random matrix whose columns
are orthogonal to the vector $\bb$. It is determined by first generating a matrix 
$\Omega_{\ell-1}$, whose entries come from a normal distribution 
with mean zero and variance one. Then, we orthonormalize the columns of
$\Omega_{\ell-1}$ with respect to 
$\bb/\|\bb\|_2$ to obtain the matrix $V_1$ in \eqref{startBlock}. In practice,
this can be achieved by applying the modified Gram-Schmidt \edit{(MGS)}
process or a QR factorization; see Algorithm~\ref{RCBA}.

It follows from the recursion relations of Algorithm~\ref{RCBA} that the block columns
$V_1,\, V_2,\, \dots,\,V_{q}$ of the matrix $\mathcal{V}_{q}$ form an orthonormal 
basis for the block Krylov subspace
\begin{equation*}
{\K}_q(A,V_1)={\rm block\,span}\lbrace V_1,AV_1,\ldots,A^{q-1}V_1 \rbrace.
\end{equation*} 
We assume throughout this paper that all the upper triangular matrices
$H_{j+1,j}$, for $j=1,2,\dots,q$, generated by
Algorithm~\ref{RCBA}, are nonsingular. This is the generic situation. 
\edit{When a subdiagonal matrix $H_{k+1,k}\in\R^{\ell\times\ell}$ is singular,
the columns of the generated matrix $V_{k+1}\in\R^{n\times\ell}$ are not 
orthogonal and 
the algorithm has to be modified. One may, for instance, reduce the
block size $\ell$ and proceed with the computations with the smaller block
size. The handling of singular blocks $H_{j+1,j}$ in the context of the
block Lanczos algorithm is discussed in \cite{Ba,baglama2003irbl}. One
may proceed analogously in the block Arnoldi algorithm. However, the 
presence of singular blocks $H_{j+1,j}$ is rare. We have not encountered
singular blocks in any of our numerical experiments and therefore will not
dwell on this issue further. Another numerical issue that we would like to 
mention is that it is straightforward to implement reorthogonalization 
in Algorithm~\ref{RCBA} by modifying line~\ref{line11}. However, we found
that reorthogonalization is not needed to improve the quality of the computed 
solutions.}

\edit{We offer some heuristic guidance for choosing the block size $\ell$.
Let the matrix $Q_\ell\in\R^{n\times\ell}$ be the same as in 
\eqref{startBlock0}. Halko et al. \cite[Corollary\,10.9]{HMT} bound
$\|A-Q_\ell Q_\ell^TA\|_2$ in terms of the singular values 
$\sigma_1\ge\sigma_2\ge\ldots\sigma_n\ge 0$ of $A$. Here, and
throughout this paper $\|\cdot\|_2$ denotes the Euclidean vector norm or
the spectral matrix norm. Specifically, choose a target rank $k\ge 2$ and an
oversampling parameter $p\ge 4$ with $\ell=k+p\le n$. Then
\[
\sigma_{\ell+1}\le\|A-Q_\ell Q_\ell^TA\|_2\leq
\left(1+8\sqrt{\ell p \log p}\right)\sigma_{k+1}+
3\left(\ell\sum_{j>k}\sigma_j^2\right)^{1/2}
\]
with probability not less than $1-3p^{-p}$. To make the upper bound small, $k$
has to be large enough so that $\sigma_{k+1}$ is small and $p$ has to be 
large enough, say $p\ge 4$, so that the bound holds with high probability. 
Moreover, Tropp and Webber \cite[Lemma\,5.1]{TW} show that for orthogonal
projectors $P_1,P_2\in\R^{n\times n}$ such that $\mathcal{R}(P_1)\subset
\mathcal{R}(P_2)$, it holds
\[
\|A-P_2A\|_r\leq\|A-P_1A\|_r
\]
for any Schatten $r$-norm with $1\le r\le\infty$. These results indicate that
the faster the singular values $\sigma_k$ of $A$ decay to zero with increasing
index $k$, the smaller $\ell$ can be chosen, and the fewer steps $q$ with the
RCBA algorithm are required to achieve an accurate approximation of $A$.}

When the matrix $A$ is symmetric, Algorithm~\ref{RCBA} simplifies to the 
Randomized Classical Block Lanczos (RCBL) algorithm. The upper block Hessenberg 
${\mathcal H}_{q}$ matrix in \eqref{RBArelation} simplifies to a block 
tridiagonal matrix with a symmetric leading principal submatrix 
${\mathcal H}_{q,q}\in\R^{\ell q\times \ell q}$. In particular, the recursion 
formulas in the RCBA algorithm simplify. The RCBL algorithm is described by 
Algorithm~\ref{CRBL}.

\begin{algorithm}[ht!] \label{CRBL}
\caption{RCBL method}
\textbf{Input:} $A\in\R^{n\times n}$, $\bb\in\R^{n}$, block size $\ell \in \mathbb{N}$, number of steps $q \in \mathbb{N}$\;
\textbf{Output:} $\mathcal{H}_{q}\in\R^{\ell (q+1)\times \ell q}$, $\mathcal{V}_{q+1}\in\R^{n\times\ell(q+1)}$ which satisfy \eqref{RBArelation}\;
    Generate \edit{$\Omega_{\ell-1} \in\R^{n\times (\ell-1)}$ with $\mathcal{N}\left(0,1\right)$ entries}\;
    Prepend $\bb$ to $\edit{\Omega_{\ell-1}}$:\quad\edit{$\Omega_{\ell}= \left[\bb,\,\Omega_{\ell-1}\right]$}\;
    Normalize: set $[V_1,\sim]=\operatorname{qr}(\edit{\Omega_{\ell}},\textrm{`econ'})$ \edit{or via MGS}\;
    Set $V_0=0$, $H_{0,1}=0$\;
    \For{$j=1,2,\ldots, q$}{
			$W =AV_j - V_{j-1}H_{j-1, j}$\;
            $H_{j,j}=V_j^TW$\;
		    $W=W-V_jH_{j,j}$\;
		    $[V_{j+1},H_{j+1,j}]=\operatorname{qr}(W,\textrm{`econ'})$\;
		    $H_{j,j+1}=H_{j+1,j}^T$\;
    	}
\end{algorithm}

\subsection{Global and loop-interchange variants of block Krylov algorithms}
\label{global_LI}
This subsection describes the randomized global and loop-interchange block 
Krylov methods. Their performance for approximating \eqref{fab} will be 
illustrated in Section~\ref{sec:4}. Non-randomized global block Krylov methods 
were first discussed by Jbilou et al.~\cite{Jbilou99,JST}, who were interested 
in solving linear systems of equations with a large matrix and several 
right-hand side vectors. \edit{In iterative block Krylov methods, global 
methods have an advantage over classical block Krylov methods that the reduced
matrix determined by the former methods may have a smaller condition number 
than the matrix of the latter methods after the same number of steps $q$. This 
is important when solving linear discrete ill-posed problems for which the
matrix $A$ is severely ill-conditioned and the right-hand side vectors are
contaminated by measurement errors; see \cite{Beik11,BEGJOR,BOR,Onisk23} for 
further discussion and applications. Frommer et al.~\cite{FLS} investigate the
application of global methods to the approximation of matrix functions of the
form $f(A)B$, where $B$ is a block vector.}

Introduce the inner product 
\[
\langle M_1,M_2\rangle_F=\operatorname{trace}(M_1^T M_2), \quad\text{for }M_1,M_2\in\R^{n\times\ell}.
\] 
We note that $\langle M_1,M_1\rangle_F=\|M_1\|_F^2$ where $\|\cdot\|_F$ denotes
the matrix Frobenius norm. 

The global block Arnoldi method is obtained by using the above inner product in
the Arnoldi iteration and require the generated block vectors 
$V_1,V_2,\ldots,V_{q+1}$ to be F-orthonormal, i.e.,
\begin{equation} \label{fOrthog}
    \left\langle V_i,V_j \right\rangle_F = \begin{cases}
        1, \quad i=j,\\
        0, \quad i \neq j.
    \end{cases}
\end{equation}
The global block Arnoldi decomposition, which is analogous to 
\eqref{RBArelation}, is given by
\begin{equation} \label{RBArelation long}
A\mathcal{V}_q = \mathcal{V}_{q+1}\mathcal{H}_{q}.
\end{equation} 
Here, $\mathcal{V}_{q+1}=[V_1,V_2,\ldots,V_{q+1}]\in\R^{n \times \ell(q+1)}$ 
and $\mathcal{V}_q$ contains the first $q$ block vector of $\mathcal{V}_{q+1}$.

We can define the initial block vector of the matrix ${\mathcal V}_q$ either by
\eqref{startBlock0} or \eqref{startBlock}; in our numerical experiments we 
focus on the latter definition. The matrix 
$\mathcal{H}_{q}\in\R^{\ell(q+1)\times \ell q}$ in \eqref{RBArelation long} is 
of block upper Hessenberg form and has block entries defined by
\begin{equation} \label{globalProduct}
H_{i,j}=\operatorname{trace}(V_i^TAV_j)I_\ell \in \mathbb{R}^{\ell \times\ell}.
\end{equation}
Hence the matrix $\mathcal{H}_q$ is sparse. The algorithm so obtained is referred to as the randomized global block Arnoldi (RGBA) method and is described by Algorithm~\ref{alg:3}.
When the matrix $A$ is symmetric, the RGBA method simplifies to the randomized global 
block Lanczos (RGBL) method. We refer to Section~\ref{sec:unified} for the algorithm.

\begin{algorithm}[ht!] \label{alg:3}
\caption{RGBA method}
\textbf{Input:} $A\in\R^{n\times n}$, $\bb\in\R^{n}$, block size $\ell \in \mathbb{N}$, and number of steps $q \in \mathbb{N}$\;
\textbf{Output:} $\mathcal{H}_{q}\in\R^{\ell (q+1)\times \ell q}$, $\mathcal{V}_{q+1}\in\R^{n\times\ell(q+1)}$ which satisfy \eqref{RBArelation long}\;
    Generate \edit{$\Omega_{\ell-1} \in\R^{n\times (\ell-1)}$ with $\mathcal{N}\left(0,1\right)$ entries}\;
    Prepend $\bb$ to $\edit{\Omega_{\ell-1}}$:\quad\edit{$\Omega_{\ell}= \left[\bb,\,\Omega_{\ell-1}\right]$}\;
    Normalize: $V_1=\edit{\Omega_{\ell}}\left(\|\edit{\Omega_{\ell}}\|_FI_\ell\right)^{-1}$\;
    \For{$k=1,\ldots,q$}{
        Compute $W=AV_k$\;
            \For{$j=1,\ldots,k$}{
                $H_{j,k}=\operatorname{trace}(V_j^TW)I_\ell$\;
                $W=W-V_jH_{j,k}$\;
            }
        $H_{k+1,k}=\Vert W\Vert_F I_\ell$\;
        $V_{k+1}=WH_{k+1,k}^{-1}$\;
    }
\end{algorithm}

The first discussion of the non-randomized loop-interchange block method is
probably due to Rashedi et al.~\cite{Rashedi16}. More recently, this block 
method has been considered by Frommer et al.~\cite{FLS}.
The loop-interchange block method enforces orthogonality between the $k$th 
columns of arbitrary block vectors $V_i$ and $V_j$ for 
$i,j \in \{1,2,\dots,q+1\}$, i.e., 
$$V_i(:,k)^TV_j(:,k)=\begin{cases}
    1,\quad i=j,\\ 0,\quad i\neq j,
\end{cases} \qquad \text{for}\,\,k=1,2,\dots,\ell.$$
Here, $V_p(:,q) \in \mathbb{R}^{n}$ denotes the $q^{th}$ 
column of the $p^{th}$ block vector of $\mathcal{V}_{q+1}$. 

The block entries $H_{i,j}\in\R^{\ell\times \ell}$ of $\mathcal{H}_q$ are
determined by ${\rm diag}\left(V_i^TV_j\right)$, where ${\rm diag}(\cdot)$ is 
the operator that sets all off-diagonal entries of a given matrix to zero, 
while preserving its diagonal. Similarly to the RGBA method, the
$q^{th}$ iterate of the randomized loop-interchange block Arnoldi decomposition
produces an analog of the decomposition \eqref{RBArelation long}, with the same
sparsity pattern. The initial block vector of the matrix ${\mathcal V}_q$ is 
given by either \eqref{startBlock0} or \eqref{startBlock}. Differently from 
the global variant, the $\ell$ diagonal entries of $H_{i,j}$ are generally not 
the same. Instead they take the form
\begin{equation} \label{loopProduct}
H_{i,j}=\operatorname{diag}(V_i^TAV_j)=\begin{bmatrix}
    V_i(:,1)^TAV_j(:,1) & & \\
    & \ddots & \\
     & & V_i(:,\ell)^TAV_j(:,\ell) \\
\end{bmatrix}.
\end{equation}
We refer to this block Arnoldi method and its
Lanczos counterpart as randomized loop-interchange block Arnoldi (RLBA) and
randomized loop-interchange block Lanczos (RLBL) methods, respectively. We 
refer to Section~\ref{sec:unified} for algorithms.

\edit{We note that in the global or loop-interchange block Krylov methods with 
block size $\ell$ a singular sub-diagonal block $H_{k+1,k}$ may occur, at least
in theory. For the global method, this breakdown happens when $\|W\|_F$ from line 12 of 
Algorithm \ref{alg:3} vanishes which indicates that the grade of $V_1$ is $\ell$; see \cite[Proposition 2]{Jbilou99}. 
In fact, one of the reasons for using the global method is that singular and nearly singular 
sub-diagonal blocks occur exceedingly seldom.
In the case of the loop-interchange method, one has to 
check whether any diagonal entry of the matrix 
$\text{diag}(\|w_1\|_2,\dots,\|w_{\ell}\|_2)$ is equal to zero; see line 23 of
Algorithm \ref{RBK_long}. Here, the $w_i$, $i=1,\dots,\ell$, denote the columns
of the matrix $W$. We did not experience any issues with singular or nearly 
singular sub-diagonal blocks in numerical experiments with the global and
loop-interchange methods. We will not discuss this situation further, but we note 
that if a computed sub-diagonal block for the loop-interchange is (nearly) singular, then it is possible 
to reduce the block size so that computation can be continued. Since, however,
it can be complicated to change the block size, we recommend that the user
instead restarts the algorithm with a smaller $\ell$.}

\subsection{A unified framework for the RBA and RBL methods}\label{sec:unified}
We present a unified framework for the randomized block Arnoldi (RBA) and 
randomized block Lanczos (RBL) methods discussed in the previous two sections. 
The randomized block classical, global, and loop-interchange methods differ in
their choices of \emph{block products} and \emph{scaling functions} that they
employ. We use the naming convention of Frommer et al.~\cite{FLS}.

Each one of the classical, global, and loop-interchange methods considered
in this paper employ a block product, which we will denote by 
$\textrm{prd}(M_1,M_2) \in \mathbb{R}^{\ell \times \ell}$ for two generic 
block vectors $M_1,\,M_2 \in \mathbb{R}^{n \times \ell}$. The block products
of the RCBA and RCBL methods, for example, are defined by 
\[
\textrm{prd}(M_1,M_2) = M_1^T M_2.
\]
The block products for the global and loop-interchange methods are given by 
\eqref{globalProduct} and \eqref{loopProduct}, respectively. They also are 
defined in Table~\ref{INP} for completeness.

A block vector, $M_1 \in \mathbb{R}^{n \times 
\ell}$ of full rank is said to be \emph{normalized} if 
\[\textrm{prd}(M_1,M_1)=I_{\ell}.\]
To normalize an arbitrary block vector $M_1$, we define the 
scaling function $N(M_1) \in \mathbb{R}^{\ell \times \ell}$. For example,
the scaling function for the RCBA and RCBL methods is given by the upper 
triangular matrix $R\in\R^{\ell\times\ell}$ of the ``skinny'' QR factorization 
$M_1=QR$, where the matrix $Q\in\R^{n\times\ell}$ has orthonormal columns.
We refer to Table~\ref{INP} for the definitions of the block vector products 
and scaling functions used for methods discussed. Algorithm~\ref{RBK_long} 
provides a general algorithm for all the block methods considered. 
\edit{We show in Section~\ref{sec:3} how the matrices $\mathcal{V}_{q+1}$ and
$\mathcal{H}_q$ are used to approximate $f(A)\bb$.}

\begin{table}[!ht]
\centering
\begin{tabular}{c|c|c}
\hline
\makecell{\textbf{Block method}\\\textbf{name}} & \makecell{\textbf{Block vector product}\\ $\text{prd}(M_1,M_2)$} & \makecell{\textbf{Scaling function}\\ $\textrm{N}(M_1)$} \\
\hline\hline
Classical & $M_1^T M_2$ & \makecell{$R$\,\,($^{\ast}$see caption)} \\
\hline
Global & $\operatorname{trace}(M_1^T M_2) I_\ell$ & $\|M_1\|_F I_{\ell}$ \\
\hline
\makecell{Loop-\\Interchange} & $\operatorname{diag}(M_1^T M_2)$ & $\operatorname{diag}\big(\|m_1\|_2, \ldots, \|m_\ell\|_2\big)$ \\
\hline
\end{tabular}
\caption{Block method names with their related block products and scaling
functions. $^{\ast}$Note that the scaling function for the classical method is
the factor $R\in \mathbb{R}^{\ell \times \ell}$ from the ``skinny'' QR 
factorization of $M_1$. \edit{It also can be computed by applying the MGS 
process to the columns of $M_1$}.}
\label{INP}
\end{table}

\begin{algorithm}[ht!] \label{RBK_long}
\caption{RBL and RBA methods with choice of prd($\cdot$,$\cdot$) and N($\cdot$)}
\textbf{Input:} $A\in\R^{n\times n}$, $\bb\in\R^{n}$, block size $\ell \in \mathbb{N}$, number of steps $q \in \mathbb{N}$, block product $\operatorname{prd}(\cdot,\cdot)$, and scaling function $\operatorname{N}(\cdot)$\;
\textbf{Output:} $\mathcal{H}_{q}\in\R^{\ell (q+1)\times \ell q}$, $\mathcal{V}_{q+1}\in\R^{n\times \ell(q+1)}$ satisfying generalized \eqref{RBArelation long}\;
    Generate \edit{$\Omega_{\ell-1} \in\R^{n\times (\ell-1)}$ with $\mathcal{N}\left(0,1\right)$ entries}\;
    Prepend $\bb$ to $\edit{\Omega_{\ell-1}}$:\quad\edit{$\Omega_{\ell}= \left[\bb,\,\Omega_{\ell-1}\right]$}\;
    Normalize: $V_1=\edit{\Omega_{\ell}}\big(N(\edit{\Omega_{\ell}})\big)^{-1}$\;
    \eIf{$A$ \textrm{is symmetric}}{
    Set $V_0=0$, $H_{0,1}=0$\;
    \For{$j=1,2,\ldots, q$}{
			$W =AV_j - V_{j-1}H_{j-1, j}$\;
            $H_{j,j}=\textrm{prd}(V_j,W)$\;
		    $W=W-V_jH_{j,j}$\;
		    $H_{j+1,j}=N(W)$\;
		    $V_{j+1}=WH_{j+1,j}^{-1}$\;
		    $H_{j,j+1}=H_{j+1,j}^T$\;
    	}\label{line15}
    	}{
		\For{$k=1,\ldots,q$}{
			Compute $W=AV_k$\;
				\For{$j=1,\ldots,k$}{
					$H_{j,k}=\operatorname{prd}(V_j,W)$\;
					$W=W-V_jH_{j,k}$\;
				}
			$H_{k+1,k}=\operatorname{N}(W)$\;
			$V_{k+1}=WH_{k+1,k}^{-1}$\;
		}\label{line25}
		}
\end{algorithm}


\section{Approximation of $f(A)\bb$ by randomized block Krylov methods} \label{sec:3}
\edit{We discuss how to efficiently compute approximations of \eqref{fab} by 
using the randomized block Krylov methods described in the previous section. 
The key idea is that if a polynomial $p$ of degree at most $q-1$ approximates 
$f$ accurately on the region of the spectrum of $A$ relevant to $\bb$, then 
$f(A)\bb \approx p(A)\bb$. Theorem~\ref{globalprf} below shows that whenever 
$p$ is such a polynomial, the vector $p(A)\bb$ can be computed exactly from the
small matrix $\mathcal{H}_{q,q}$ and the basis $\mathcal{V}_q$ generated by 
Algorithm~\ref{RBK_long} without explicitly forming $p(A)$ or the powers 
$A^j\bb$. This motivates the approximation}
\[
f(A)\bb \approx \mathcal{V}_q f(\mathcal{H}_{q,q}) \be_1 \xi,
\]
\edit{which we use throughout the remainder of this section and in the computed examples of Section~\ref{sec:4}. The scalar $\xi$ is a scaling factor. The 
precise statement, given below, is a block generalization of Saad's work 
\cite{saad1992_FOM}; see also \cite{BR}.
} 

\begin{theorem}\label{globalprf}
Let $\mathcal{V}_q=[V_1,V_2, \ldots,V_q]$ and $\mathcal{H}_{q,q}$ be the 
result of an application of $q$ steps of Algorithm~\ref{RBK_long} to the 
matrix $A\in\R^{n\times n}$ with initial vector $\bb\in\R^{n}$. Let $\xi>0$ 
be a suitable scalar determined by the randomized block Krylov method (see 
the end of this section for a discussion on $\xi$). Then for any monomial 
$p_j$ of degree $0\le j\le q-1$, we have
\begin{equation} \label{proj_fab}
p_j(A)\bb=\edit{\mathcal{V}_qp_j(\mathcal{H}_{q,q})} \be_1\xi,
\end{equation}
where $\mathcal{V}_q$ and $\mathcal{H}_{q,q}$ satisfy the relation 
\[
A\mathcal{V}_q=\mathcal{V}_q\mathcal{H}_{q,q}+
V_{q+1}H_{q+1,q}\mathbf{E}_q^T,
\]
which is analogous to \eqref{RBArelation} or \eqref{RBArelation long}.
Here $\mathbf{E}_q \in \mathbb{R}^{\ell q \times \ell}$ is made up of the 
last $\ell$ columns of $I_{\ell q} \in \mathbb{R}^{\ell q \times \ell q}$.
\end{theorem} 

\noindent\text{Proof.} We proceed by induction over the degree of $p_j$. For 
$j=0$, we have
\[
p_{0}(A)\bb=A^{0}\bb=\frac{\bb}{\xi}\xi =\edit{\mathcal{V}_q} \be_1\xi =\edit{\mathcal{V}_qp_0(\mathcal{H}_{q,q})}\be_1\xi.
\]
For the inductive step, assume that 
$p_j(A)\bb=\edit{\mathcal{V}_qp_j(\mathcal{H}_{q,q})} \be_1\xi$ \edit{holds 
for all polynomials of degree $j$} and consider 
\begin{align*}
p_{j+1}(A)\bb&=AA^{j}\bb\\
&=A\edit{\mathcal{V}_qp_j(\mathcal{H}_{q,q})}  \be_1\xi\\
&=(\edit{\mathcal{V}_q\mathcal{H}_{q,q}+V_{q+1}H_{q+1,q}\mathbf{E}_q^T})p_j(\edit{\mathcal{H}_{q,q}})\be_1\xi\\
&=\edit{\mathcal{V}_q\mathcal{H}_{q,q}}p_j(\edit{\mathcal{H}_{q,q}})\be_1\xi+\edit{V_{q+1}H_{q+1,q}\mathbf{E}_q^T}p_j(\edit{\mathcal{H}_{q,q}})\be_1\xi.
\end{align*}
Since $\mathbf{E}_{q}^Tp_j(\mathcal{H}_{q,q}) \be_1={\bf 0}$, the second term 
vanishes. This gives 
\[
p_{j+1}(A)\bb=\edit{\mathcal{V}_qp_{j+1}(\mathcal{H}_{q,q})}\be_1\xi.
\]
We conclude that the relation \eqref{proj_fab} holds for every integer
$0\le j\leq q-1$.  $\hfil\Box$

We turn to the choice of the scaling factor $\xi>0$ in \eqref{proj_fab}. 
Since the first column of the matrix $\mathcal{V}_q$ for any of the methods 
discussed is a scaling of $\bb$, we need choose the scalar $\xi$ so that 
the first column of $\mathcal{V}_q\xi$ equals $\bb$. For the classical and 
loop-interchange methods, we let $\xi=\Vert \bb\Vert_2$ and for the global
methods we use $\xi=\Vert V_1\Vert_F$.

\section{Computed examples}\label{sec:4}
This section presents numerical results for the randomized block Arnoldi and 
Lanczos methods. \edit{We chose the functions $f$ in this section to illustrate
the performance of our methods across a representative range of applications, 
namely regularization, inverse problems, a square-root function, and network 
problems. Each function can be computed directly in MATLAB (via \texttt{expm},
\texttt{sqrtm}, or backslash) for comparison purposes.} We illustrate how the
block size $\ell$ influences both the relative difference between iterations 
and the computation time (in seconds) for four representative examples:
\begin{itemize}
\item a regularization problem with an ill-conditioned matrix,
\item a square-root function example with several choices of $\bb$,
\item an inverse function example, and
\item a consistency study using five-number summaries for two network problems.
\end{itemize}

To stop the iteration,
we measure the \emph{relative difference between iterations} (RDBI) 
\[
\textrm{RDBI}=\frac{\Vert f(\mathcal{H}_{q,q})\be_1 \xi - f(\mathcal{H}_{q-1,q-1}) \be_1 \xi 
\Vert_2}{\Vert f(\mathcal{H}_{q,q}) \be_1 \xi \Vert_2},
\]
where $f(\mathcal{H}_{q,q}) \be_1 \xi$ is the approximation defined in 
Theorem~\ref{globalprf} and $f(\mathcal{H}_{q-1,q-1}) \be_1 \xi$ is appended 
with zeros to make the two vectors of the same length. The approximate solution
in the original space is given by $\mathcal{V}_qf(\mathcal{H}_{q,q})\be_1\xi$.
RDBI is calculated after every depth step $q>2$ in Algorithm~\ref{RBK_long} 
(after lines~\ref{line15} or~\ref{line25} for RBL and RBA, respectively). 

For the purposes of this paper, we fix the tolerance for RDBI at $10^{-15}$
except for the regularization problems where noise is introduced. For the
latter problems, we raise the tolerance to $10^{-10}$. The methods perform
well also for larger tolerances such as $10^{-5}$ on each of the problem types. Once the chosen tolerance is reached, the algorithm is terminated.  

To account for the inherent randomness of the block randomized Krylov methods,
the first three examples were run independently 15 times and we report the 
average RBDI and timings. All computations were performed in MATLAB 2025a on a
MacBook Pro equipped with a 2.3 GHz 8-core Intel Core i9 processor and 16 GB of
RAM. The computations were carried in double precision, i.e., to about 15 
significant decimal digits.

\edit{We do not report peak memory usage in these experiments, because memory
use is dominated by storage of the matrix $A$, while the block Krylov basis 
$\mathcal{V}_{q+1}\in\R^{n\times\ell(q+1)}$ requires comparatively little 
additional memory since $\ell(q+1)\ll n$ in all cases considered. A detailed 
investigation of memory movement costs for block Krylov methods is left to 
future work.}

\subsection{Matrices that stem from the discretization of linear ill-posed 
problems}\label{Reg}  
We apply randomized block Lanczos methods to the approximation of expressions \eqref{fab}
with matrices $A$ and vectors $\bb$ that arise when solving linear
ill-posed problems.
Such problems usually stem from the discretization of
Fredholm integral equations of the first kind. Specifically, we 
consider the integral equations described by Baart~\cite{Baart}, Phillips~\cite{Ph}, and 
Shaw~\cite{Shaw}. 

The integral equation in~\cite{Baart} is discretized by a Galerkin method with piece-wise 
constant test and trial functions using software provided by Hansen~\cite{Tools}. The 
software yields a linear system of equations
\begin{equation}\label{linsys}
M\bx=\bc,
\end{equation}
with a matrix $M\in\R^{128^2\times 128^2}$ and right-hand side vector $\bc\in\R^{128^2}$.
The matrix $M$ is very ill-conditioned and in applications the vector $\bc$ is generally
contaminated by error. To avoid obtaining a large propagated error in the computed
(approximate) solution of \eqref{linsys}, Tikhonov regularization is employed. This
regularization method replaces the least-squares solution of \eqref{linsys} by
the solution of a regularized system
\begin{equation}\label{tik}
\left(M^TM+\frac{1}{\mu} I\right)\bx = M^T\bc,
\end{equation}
where $\mu>0$ is a regularization parameter. It is chosen to give a desirable
approximate solution in the presence of errors in the data vector $\bc$
and round-off errors introduced during the computations; see, e.g.,~\cite{CG,
Ha1} for discussions on the solution of ill-posed problems that arise when 
solving integral equations of the first kind and the choice of the 
regularization parameter $\mu$.

As an alternative to solving the Tikhonov problem \eqref{tik}, it was suggested in~\cite{CR02,CDD} to instead solve the matrix exponential problem 
\[
\exp(\mu M^TM){\bx}=\mu M^T\bc.
\]
The relation between this equation and \eqref{tik} can be seen by using the first two 
terms in the Maclaurin expansion of $\exp(t)$. We apply the methods of the present 
paper to evaluate \eqref{fab} with $f=\exp(-t)$, $A=\mu M^TM$, and $\bb=M^T\bc$. In the computed
examples, we set $\mu=1$. We remark that the regularization method discussed in~\cite{CR02} focuses on the situation when the matrix $M$ is symmetric. Additionally, in each of the problems in this section, we add 1\% of Gaussian white noise to the vector $\bc$ before forming $\bb=M^T\bc$.

\begin{figure}[ht]
\begin{center}
\includegraphics[width=\linewidth, height=.75\linewidth]{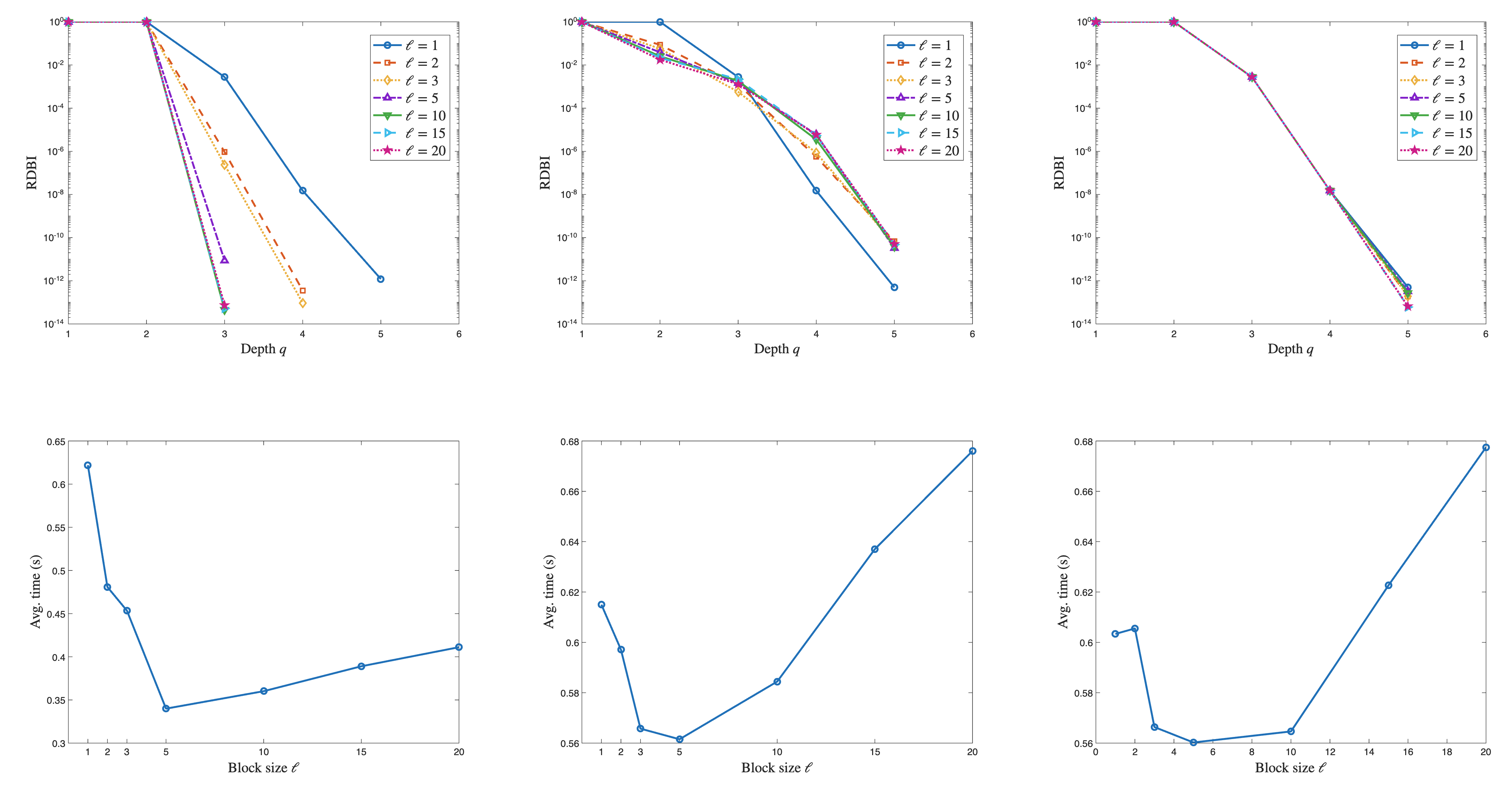}
\caption{Regularization example: RDBI for the RCBL, RGBL, RLBL methods (top row, left to 
right) and timings in seconds (bottom row) for the test problem Baart.}\label{RegBaart}
\end{center}
\end{figure}

Computed results presented in Figure~\ref{RegBaart} indicate that all RBL methods satisfy the
prescribed RDBI tolerance. Among them, RCBL achieves the solution faster than
RGBL and RLBL. Notably, RCBL, with the exception of block size $\ell=1$, outperforms the best times produced by  RGBL and RLBL, while increasing the block size to $\ell=3$ gives a 20\% reduction in runtime for RCBL over RGBL and similarly outperforms RLBL with a 18\% reduction in runtime.

\begin{figure}[ht]
\begin{center}
\includegraphics[width=\linewidth, height=.75\linewidth]{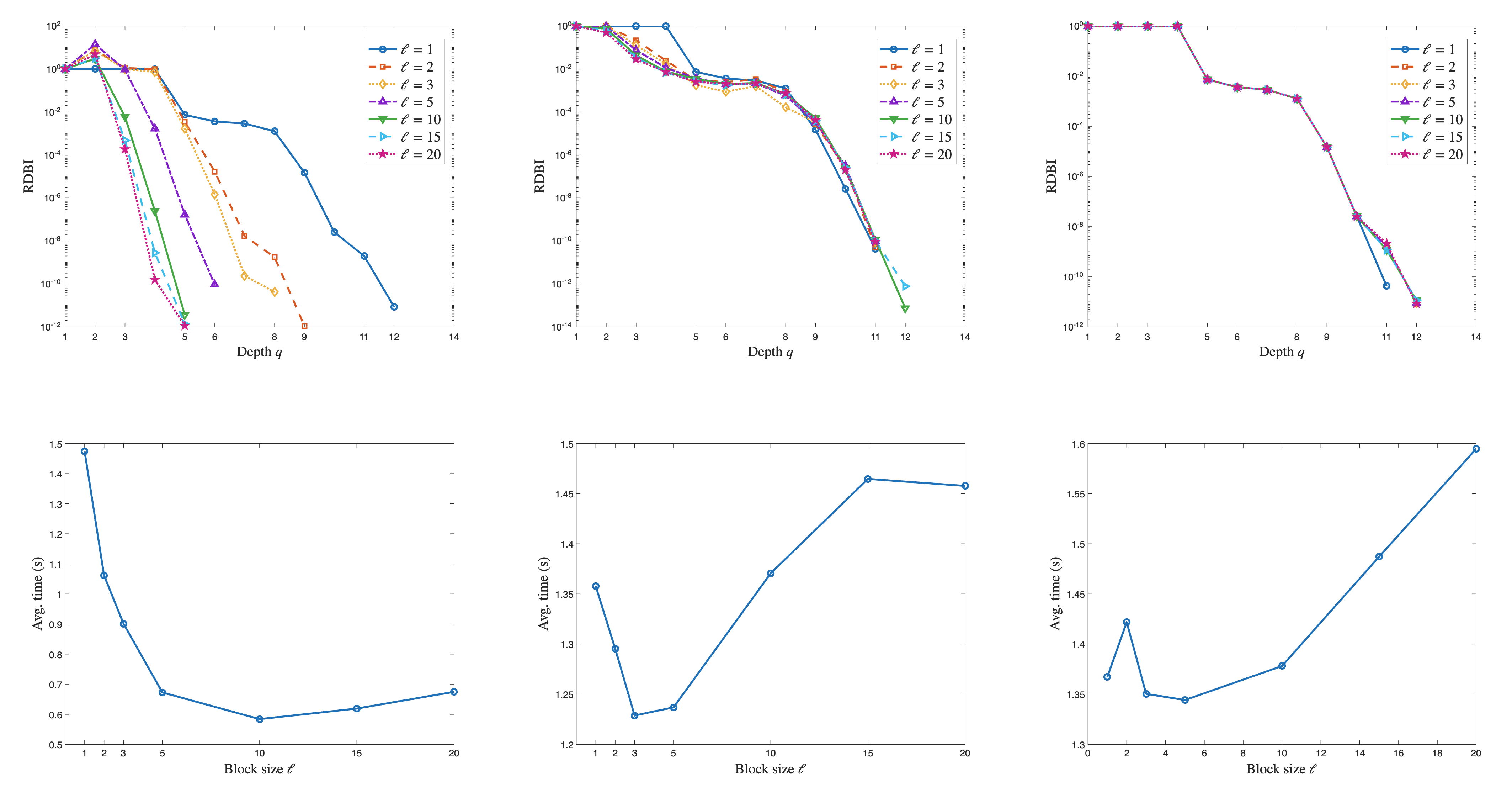}
\caption{Regularization example: RDBI for the RCBL, RGBL, RLBL methods (top row, left to
right) and timings in seconds (bottom row) for the test problem Phillips.}
\label{RegPhillips}
\end{center}
\end{figure}

Similarly, the integral equation discussed by Phillips~\cite{Ph} may be 
discretized by a Galerkin method with orthonormal box functions to obtain a 
linear system of equations \eqref{linsys} with a matrix of size 
$128^2\times 128^2$. Code for determining the matrix 
and right-hand side can be found in~\cite{Tools}. Figure~\ref{RegPhillips}
shows that, similarly to the Baart test problem, the RDBI for each RBL method
meets our tolerance. This time we see that RCBL benefits from increasing block
size, with block size 10 offering a 61\% decrease in runtime with respect to
block size 1.

\begin{figure}
\begin{center}
\includegraphics[width=\linewidth, height=.75\linewidth]{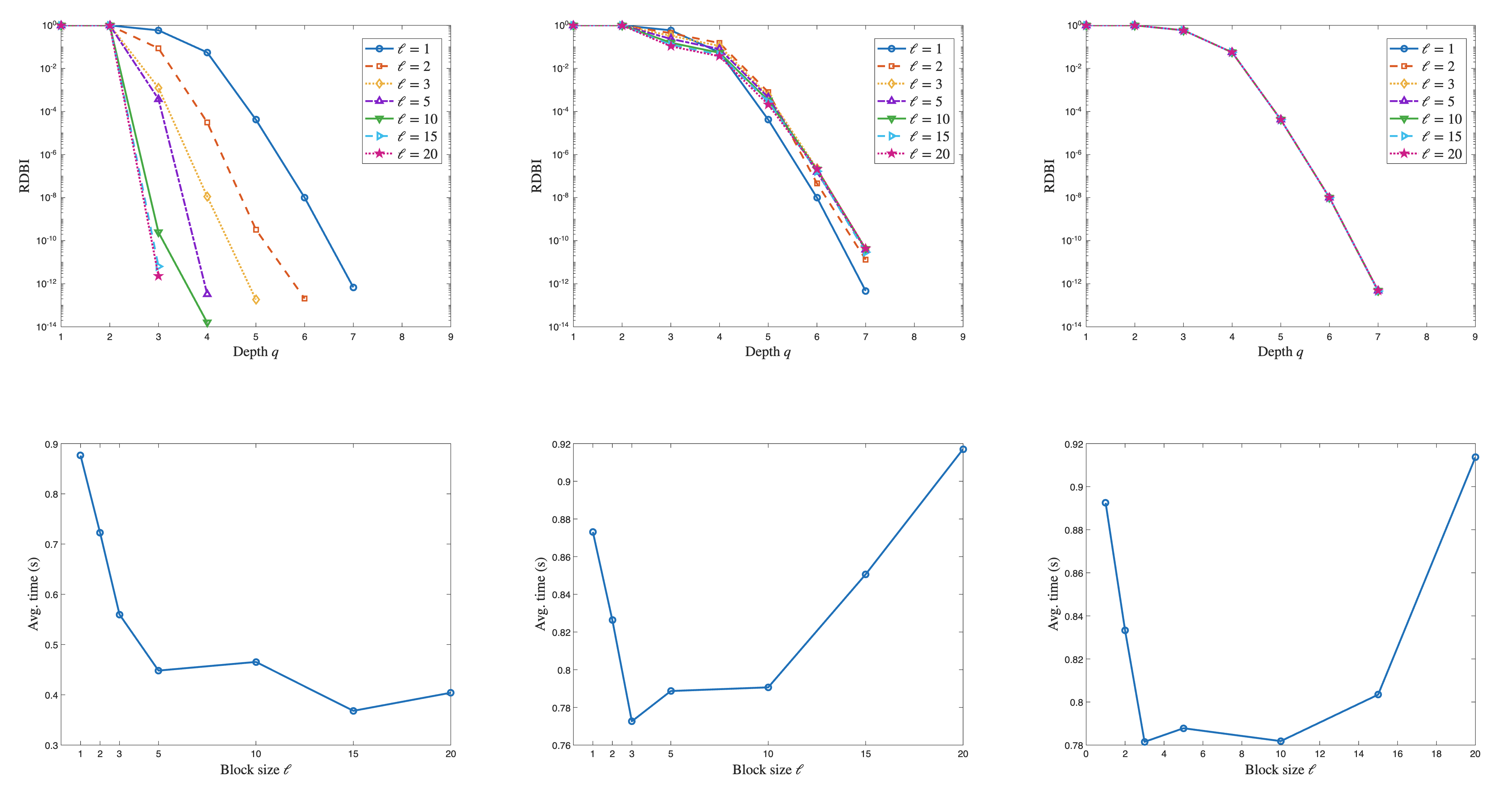}
\caption{Regularization example: RDBI for the RCBL, RGBL, RLBL methods (top row, left to
right) and timings in seconds (bottom row) for the test problem Shaw.}\label{RegShaw}
\end{center}
\end{figure}

Finally, we turn to the integral equation described by Shaw~\cite{Shaw}. It is 
discretized by a simple quadrature rule using code from~\cite{Tools}. The 
system matrix is of the same size as the Phillips and Baart problems. Computed
results displayed in Figure~\ref{RegShaw} show that, similarly to the previous 
examples, each of the RBL methods satisfies the tolerance set for RDBI. Here we
see that block size 2 of the RCBL is faster than the other methods. Further, we
see that increasing the block size from 2 to 15 reduces the runtime of RCBL by
further 32\%. Overall, block size 15 provides a 56\% decrease in runtime from 
block size 1 for RCBL.

Our experiments show that increasing the block size $\ell$ does not consistently reduce 
the average RDBI for the RCBL, RGBL, and RLBL methods. However, larger block sizes often 
decrease runtime and reduce the number of steps $q$ required to achieve the target 
tolerance. In Figure~\ref{RegBaart}, all methods benefit from increasing $\ell$. RCBL 
achieves a 45\% reduction in average runtime and fewer average depth steps when $\ell$ 
increases from 1 to 5, while RGBL and LRBL each show about a 10\% average runtime 
reduction. In Figure~\ref{RegPhillips}, RGBL shows around 15\% average runtime reduction when using
a larger block size $\ell$. However, these gains are quickly lost if increasing the block size too much. For RCBL, increasing $\ell$ from 1 to 20 yields a 60\% 
reduction in average runtime and a 50\% reduction in steps $q$. For RLBL, increasing
$\ell$ from 1 to 3 produces a 8\% runtime reduction and fewer depth steps. 
In Figure~\ref{RegShaw}, RCBL shows a 58\% reduction in average runtime and fewer average depth 
steps when increasing the block size $\ell$ from 1 to 15, while both RGBL and RLBL show roughly a 10\% 
reduction in average runtime when increasing the block size from 1 to 3.

Across all experiments, RCBL consistently outperforms the other methods. While larger
block sizes often reduce runtime, the gains diminish as $\ell$ increases, and beyond a
problem-dependent threshold, further increases in $\ell$ may lead to higher runtime;
see Figure~\ref{RegBaart}.


\subsection{Square Root}\label{SQRT}
In this example, we apply RBA methods to the non-symmetric Baart problem. We set 
\begin{equation}\label{fcnbaart}
f(A) = \sqrt{A + 0.01I} 
\end{equation}
and choose $\bb$ as a normalized vector of ones. The shift ensures that all singular 
values are positive. The matrix $A\in\R^{10000\times 10000}$ is generated with code accompanying~\cite{Tools}.

\begin{figure}[ht]
\begin{center}
\includegraphics[width=\linewidth, height=.75\linewidth]{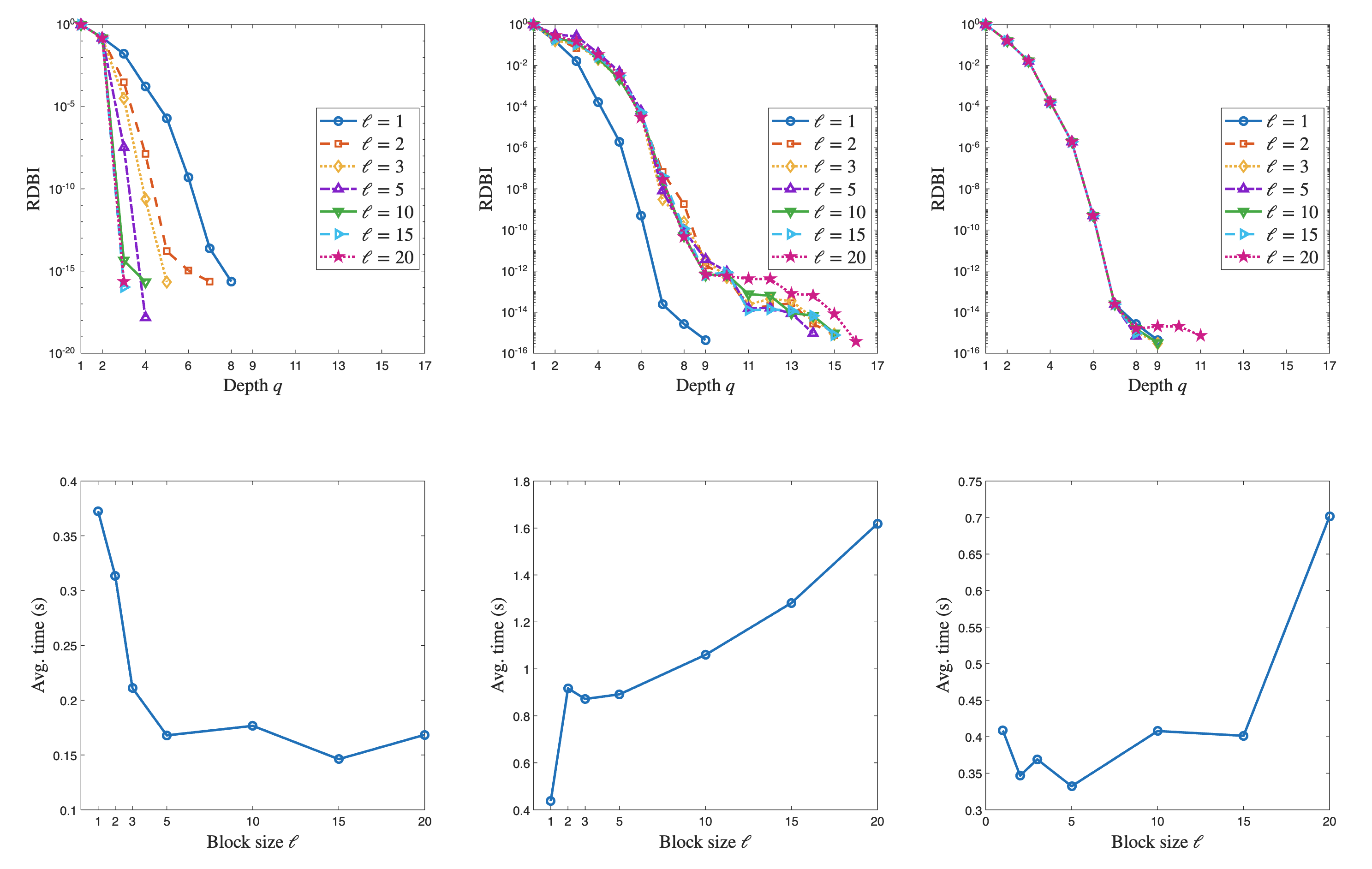}
\caption{Square Root example: RDBI for the RCBA, RGBA, and RLBA methods (top row, left to
right) and timings in seconds (bottom row) for the function \eqref{fcnbaart} with 
$\bb$ the normalized vector of all ones.}\label{SqrtBaart}
\end{center}	
\end{figure}  

Figure~\ref{SqrtBaart} shows that RBA achieves RDBI on the order of $10^{-15}$ on average.
Increasing the block size from $\ell=1$ to $\ell=15$ yields a 60\% reduction in average 
runtime and a 50\% reduction in average steps $q$ when using the RCBA. For RGBA, 
increasing the block size beyond $\ell=1$ leads to increased runtime. Finally, for the 
RLBA, increasing the block size from $\ell=1$ to $\ell=3$ yields a 17\% reduction in 
average runtime and a 15\% reduction in average steps $q$.

We note that for this example, the classical block method ultimately
outperforms the other methods in terms of runtime and depth steps $q$. 
In fact, the classical method with $\ell=1$ outperforms the best cases 
of each of the other methods.

\begin{figure}[ht]
\begin{center}
\includegraphics[width=\linewidth, height=.75\linewidth]{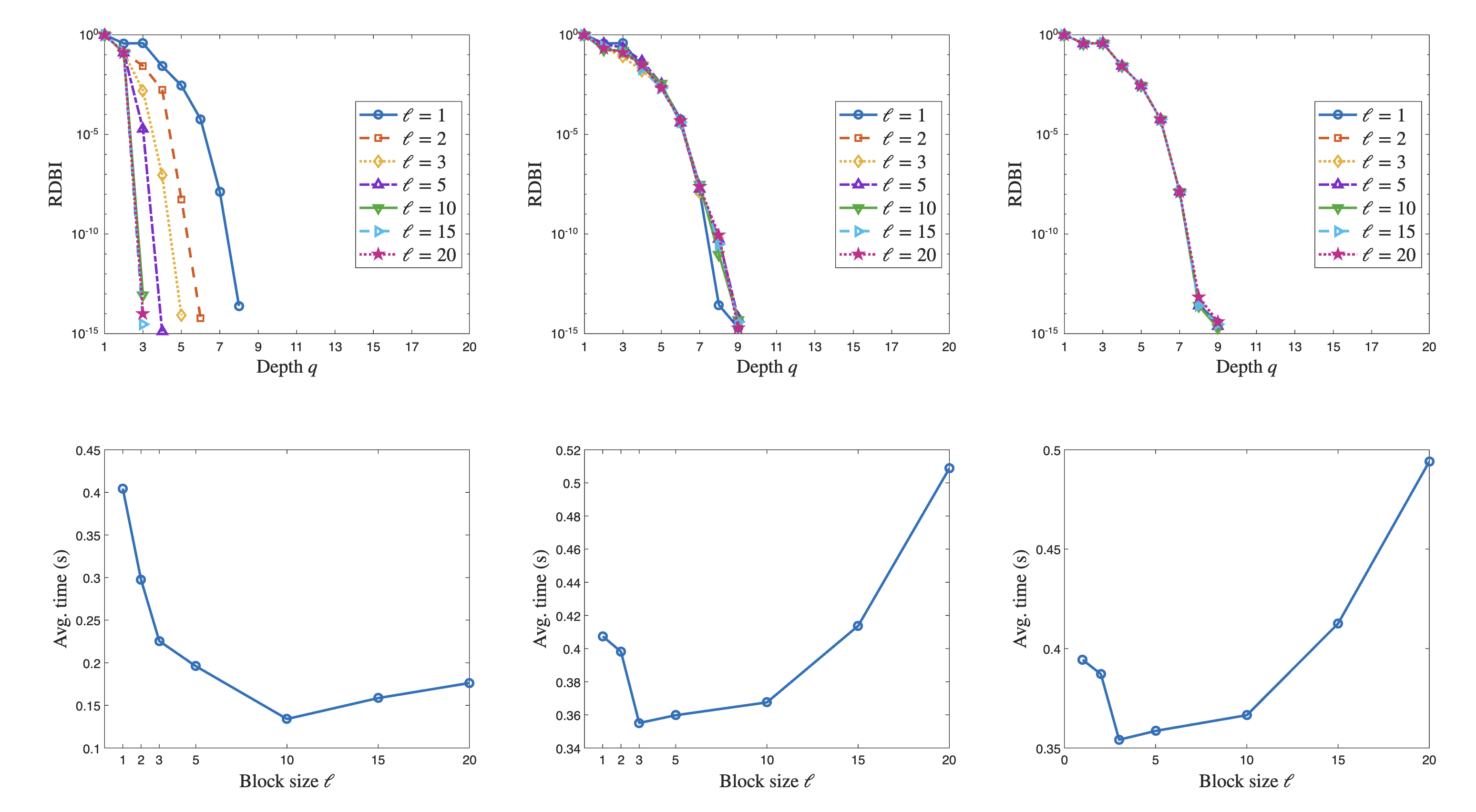}
\caption{Square Root example: RDBI for the RCBA, RGBA, and RLBA methods (top row, left to 
right) and timings in seconds (bottom row) for the function \eqref{fcnbaart} with $\bb$
a random vector.}\label{SqrtBaart2}
\end{center}
\end{figure}

In Figure~\ref{SqrtBaart2}, we repeat the experiment above with the same matrix, but with
$\bb$ chosen as a random vector whose entries are drawn from $\mathcal{N}(0,1)$. Here, we found
that RCBA benefits from increasing the block size from $\ell = 1$ to $\ell = 5$, yielding
less than a 49\% reduction in average runtime and a 44\% reduction in steps $q$. In 
contrast, increasing $\ell$ does not decrease runtime for RLBA or RGBA, though both see 
roughly a 10\% reduction in depth steps $q$. The RCBA method with $\ell=3$ outperforms the
other methods, though RLBA with $\ell=5$ performs comparably on average.

These results suggest that \edit{for RCBA} increasing the block size often reduces both runtime and depth 
steps $q$ in problems involving a square-root. \edit{However, this effect is not present in the RGBA and RLBA methods.} Additionally, the relative performance of the randomized block methods considered not only depends on the matrix $A$, but also on the vector $\bb$.


\subsection{Inverse}\label{INV}  
We next apply the RBL methods to the symmetric Shaw problem with $f(A) = (A + 0.01I)^{-1}$ and 
$\bb$ which is drawn from a standard normal distribution. The shift ensures that all eigenvalues are bounded away from zero. We let $A\in\R^{10000 \times 10000}$.

\begin{figure}[!ht]
\begin{center}
\includegraphics[width=\linewidth, height=.75\linewidth]{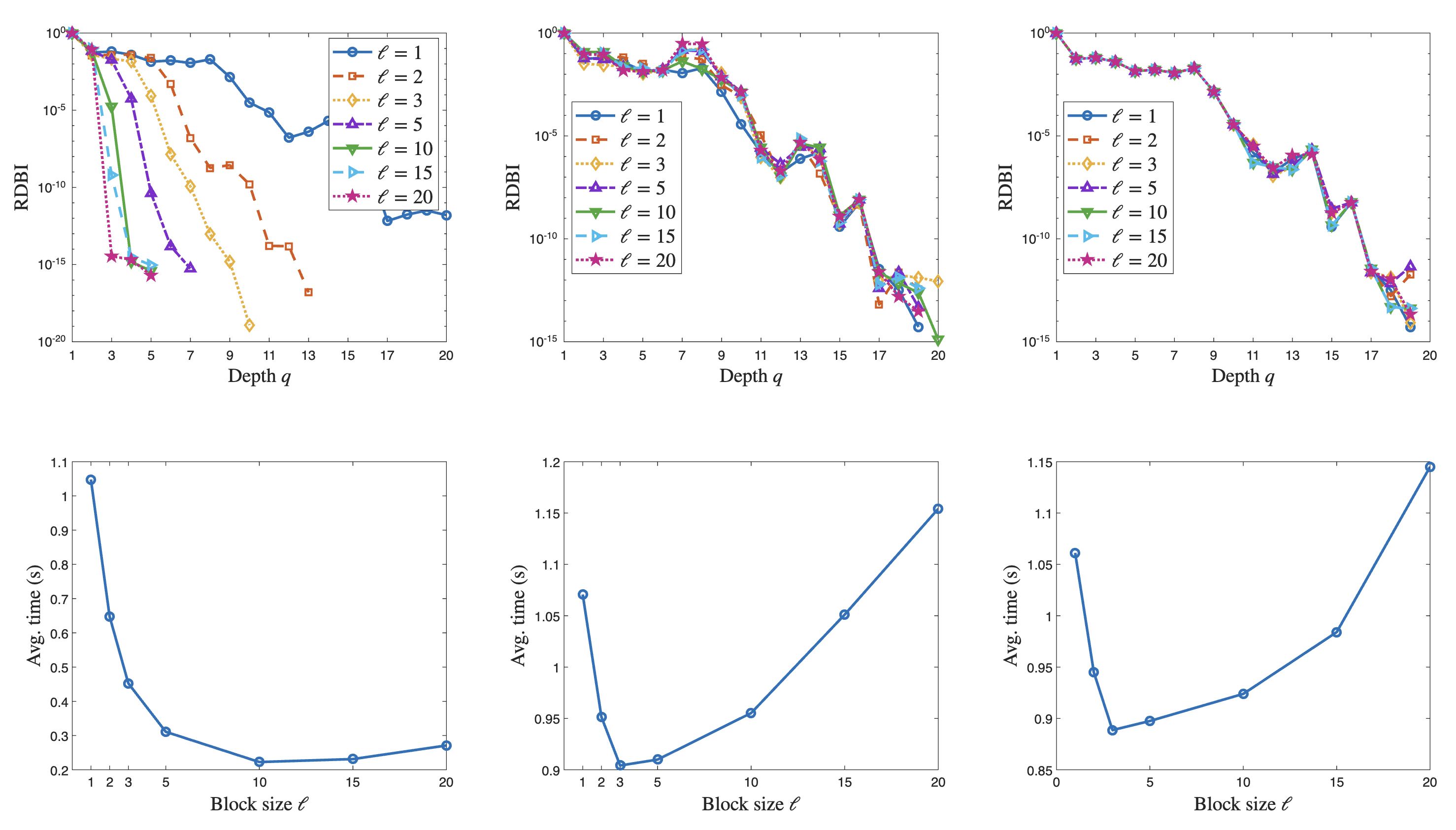}
\caption{Inverse example: RDBI for the RCBL, RGBL, RLBL methods (top row, left to right) 
and timings in seconds (bottom row) for the symmetric test problem Shaw}\label{Invfig}
\end{center}
\end{figure}  

We see in Figure \ref{Invfig}, similarly to the previous examples, that 
increasing the block size yields substantial improvements for 
RCBL: increasing $\ell$ from 1 to 15 results in a 70\% reduction in runtime and steps $q$.
RGBL shows no improvement with larger $\ell$, while LRBL gains a modest 15\% runtime 
reduction with no change in steps $q$. Overall, RCBL outperforms the other methods, providing substantial reductions in runtime and
depth steps.


\subsection{Consistency of RBL methods}
Finally, we examine the consistency of the RBL and RBA methods. We first present 
five-number summaries for 15 independent runs where $A\in\R^{22963\times 22963}$ is the symmetric adjacency matrix that is associated with the network \textit{Internet}; see~\cite{FRR22} for a discussion of the problem and associated code. The block size is fixed to be $\ell = 3$. For this problem the vector $\bb$ has all entries equal and positive and the $\ell-1$ random vectors of the initial block $\Omega_{\ell - 1}$ are also regenerated for each run.

\begin{table}[!ht]
\centering
\caption{Consistency Example: Five-number summaries over 15 runs of RBL methods for block
size $\ell=3$ applied to the symmetric network Internet problem. Time shows the execution
time of the corresponding method in seconds.}
\label{Int5Num}
\begin{tabular}{l|c|c|c|c|c|c}
\textbf{Method}&\textbf{Value}&\textbf{Min}&\textbf{Median}&\textbf{Max}&\textbf{Mean}& 
\textbf{SD}\\\hline\hline
      RCBL           &RDBI         &2.3768e-16  &5.1419e-16     &9.7191e-16  &5.5961e-16   &2.0916e-16\\
                     &Time          &0.0583      &0.0805         &0.1536      &0.0846       &0.0215\\\hline
      RGBL           &RDBI         &3.1519e-16  &7.0370e-16     &9.7148e-16  &6.9508e-16   &1.9217e-16\\
                     &Time          &0.0785      &0.0878         &0.0956      &0.0871       &0.0049\\\hline
      RLBL           &RDBI         &4.0449e-16  &7.4674e-16     &9.5794e-16  &7.1119e-16   &1.7386e-16\\
                     &Time          &0.0984      &0.1029         &0.1249      &0.1036       &0.0104\\\hline\hline
    \end{tabular}
\end{table}

Table~\ref{Int5Num} shows that each of the RBL methods performs consistently, as indicated
by the small standard deviations of the RDBI across multiple runs. For $\ell=3$, all 
methods achieve RDBIs of similar magnitude. Although RCBL exhibits the largest variability in runtime, it also attains the smallest
observed execution time. Moreover, its median runtime indicates that RCBL is typically the
fastest method, while RGBL and RLBL are slightly slower but with more consistent runtimes.

\begin{table}[!ht]
\centering
\caption{Consistency Example: Five-number summaries over 15 runs of RBA methods for block 
size $\ell=3$ applied to the non-symmetric network Wiki-votes. Time is in seconds.}\label{Int5Num2}
    \begin{tabular}{l|c|c|c|c|c|c}
      \textbf{Method}&\textbf{Value}&\textbf{Min}&\textbf{Median}&\textbf{Max}&\textbf{Mean}& \textbf{SD}\\\hline\hline
      RCBA           &RDBI         &1.5930e-16  &7.6168e-16     &2.1435e-14  &9.8521e-15   &3.4132e-14\\
                     &Time          &0.0889      &0.0914         &0.2335      &0.1017       &0.0366\\\hline
      RGBA           &RDBI         &9.3037e-16  &1.0832e-15     &1.4283e-14  &1.1286e-15   &1.5507e-14\\
                     &Time          &0.0903      &0.1143         &0.1285      &0.1128       &0.0089\\\hline
      RLBA           &RDBI         &8.5774e-16  &2.9918e-15     &7.8859e-15  &3.4337e-15   &2.1792e-15\\
                     &Time          &0.0846      &0.0897         &0.1192      &0.0924       &0.0093\\\hline\hline
    \end{tabular}
\end{table}

We repeated the experiment applying the RBA methods to the non
symmetric adjacency matrix, $A\in\R^{7115\times7115},$ associated with the undirected network
Wiki-votes; see~\cite{FRR22}.

Table~\ref{Int5Num2} highlights the various RBA methods. We note that all methods achieve
similar RDBI values and have small standard deviations for both the runtime and RDBI. These 
numbers show that the RBA method can be expected to perform well in every run. \edit{One major 
difference between Tables~\ref{Int5Num} and~\ref{Int5Num2} is that the RLBA method has 
the best time in all categories and thus would be the preferred method for this problem}

\section{Conclusion}\label{sec:5}
We considered the application of randomized block Krylov subspace
methods for the approximation of expressions of the form $f(A)\bb$.
Specifically, we considered classical, global, and loop-interchange block
methods, whose initial columns were drawn from a standard normal distribution 
and were prepended with the vector $\bb$. In our numerical experiments we found
that a value of the block size $\ell$ slightly larger than $1$ often reduces both runtime and the number of steps $q$ needed to reach a desired tolerance. \edit{However, these methods often show increased runtime as the block size grows. Heuristically, this stems from how quickly each method estimates the dominant eigenvalues values of the initial matrix $A$: the classical method estimates this information more quickly than the global or loop-interchange methods, which accounts for the runtime differences we observe in our results.} These findings underscore the practical value of randomized block Krylov methods for matrix function problems, highlighting how block size can be tuned to balance efficiency and computational cost.

\section*{Acknowledgment}
L.O. acknowledges partial support by the U.S. National Science Foundation, under grant DMS-2038118. Any opinions, findings, and conclusions or recommendations expressed in this material are those of the author(s) and do not necessarily reflect the views of the National Science Foundation. 
G. R. is partially supported by the INdAM-GNCS 2026 project ``Metodi numerici
per modelli integrali e dinamiche con memoria'' and by the research project
Fondazione di Sardegna 2024-2025 ``Integral and Discrete Inverse Problems
(InDIP)''.

\bibliographystyle{siam}
\bibliography{refs}

\end{document}